# On non-asymptotic bounds for estimation in generalized linear models with highly correlated design

Sara A. van de Geer[1]

*ETH Zürich*

**Abstract:** We study a high-dimensional generalized linear model and penalized empirical risk minimization with $\ell_1$ penalty. Our aim is to provide a non-trivial illustration that non-asymptotic bounds for the estimator can be obtained without relying on the chaining technique and/or the peeling device.

## 1. Introduction

We study an increment bound for the empirical process, indexed by linear combinations of highly correlated base functions. We use direct arguments, instead of the chaining technique. We moreover obtain bounds for an M-estimation problem inserting a convexity argument instead of the peeling device. Combining the two results leads to non-asymptotic bounds with explicit constants.

Let us motivate our approach. In M-estimation, some empirical average indexed by a parameter is minimized. It is often also called empirical risk minimization. To study the theoretical properties of the thus obtained estimator, the theory of empirical processes has been a successful tool. Indeed, empirical process theory studies the convergence of averages to expectations, uniformly over some parameter set. Some of the techniques involved are the *chaining technique* (see e.g. [13]), in order to relate increments of the empirical process to the entropy of parameter space, and the *peeling device* (a terminology from [10]) which goes back to [1], which allows one to handle weighted empirical processes. Also the *concentration inequalities* (see e.g. [9]), which consider the concentration of the supremum of the empirical process around its mean, are extremely useful in M-estimation problems.

A more recent trend is to derive non-asymptotic bounds for M-estimators. The papers [6] and [4] provide concentration inequalities with economical constants. This leads to good non-asymptotic bounds in certain cases [7]. Generally however, both the chaining technique and the peeling device may lead to large constants in the bounds. For an example, see the remark following (5).

Our aim in this paper is simply to avoid the chaining technique and the peeling device. Our results should primarily be seen as non-trivial illustration that both techniques may be dispensable, leaving possible improvements for future research. In particular, we will at this stage not try to optimize the constants, i.e. we will make some arbitrary choices. Moreover, as we shall see, our bound for the increment

---

[1]Seminar für Statistik, ETH Zürich, LEO D11, 8092 Zürich, Switzerland, e-mail: geer@stat.math.ethz.ch

involves an additional log-factor, $\log m$, where $m$ is the number of base functions (see below).

The M-estimation problem we consider is for a high-dimensional generalized linear model. Let $Y \in \mathcal{Y} \subset \mathbf{R}$ be a real-valued (response) variable and $x$ be a covariate with values in some space $\mathcal{X}$. Let

$$\left\{ f_\theta(\cdot) := \sum_{k=1}^m \theta_k \psi_k(\cdot), \theta \in \Theta \right\}$$

be a (subset of a) linear space of functions on $\mathcal{X}$. We let $\Theta$ be a convex subset of $\mathbf{R}^m$, possibly $\Theta = \mathbf{R}^m$. The functions $\{\psi_k\}_{k=1}^m$ form a given system of real-valued base functions on $\mathcal{X}$. The number of base functions, $m$, is allowed to be large. However, we do have the situation $m \leq n$ in mind (as we will consider the case of fixed design).

Let $\gamma_f : \mathcal{X} \times \mathcal{Y} \to \mathbf{R}$ be some loss function, and let $\{(x_i, Y_i)\}_{i=1}^n$ be observations in $\mathcal{X} \times \mathcal{Y}$. We consider the estimator with $\ell_1$ penalty

$$(1) \qquad \hat{\theta}_n := \arg\min_{\theta \in \Theta} \left\{ \frac{1}{n} \sum_{i=1}^n \gamma_{f_\theta}(x_i, Y_i) + \lambda_n^{\frac{2}{2-s}} I^{\frac{2(1-s)}{2-s}}(\theta) \right\},$$

where

$$(2) \qquad I(\theta) := \sum_{k=1}^m |\theta_k|$$

denotes the $\ell_1$ norm of the vector $\theta \in \mathbf{R}^m$. The smoothing parameter $\lambda_n$ controls the amount of complexity regularization, and the parameter $s$ ($0 < s \leq 1$) is governed by the choice of the base functions (see Assumption B below). Note that for a properly chosen constant $C$ depending on $\lambda_n$ and $s$, we have for any $I > 0$,

$$\lambda_n^{\frac{2}{2-s}} I^{\frac{2(1-s)}{2-s}} = \min_\lambda \left( \lambda I + \frac{C}{\lambda^{\frac{2(1-s)}{s}}} \right).$$

In other words, the penalty $\lambda_n^{\frac{2}{2-s}} I^{\frac{2(1-s)}{2-s}}(\theta)$ can be seen as the usual Lasso penalty $\lambda I(\theta)$ with an additional penalty on $\lambda$. The choice of the latter is such that adaption to small values of $I(\theta_n^*)$ is achieved. Here, $\theta_n^*$ is the target, defined in (3) below.

The loss function $\gamma_f$ is assumed to be convex and Lipschitz (see Assumption L below). Examples are the loss functions used in quantile regression, logistic regression, etc. The quadratic loss function $\gamma_f(x,y) = (y - f(x))^2$ can be studied as well without additional technical problems. The bounds then depend on the tail behavior of the errors.

The covariates $x_1, \ldots, x_n$ are assumed to be fixed, i.e., we consider the case of fixed design. For $\gamma : \mathcal{X} \times \mathcal{Y} \to \mathbf{R}$, use the notation

$$P\gamma := \frac{1}{n} \sum_{i=1}^n \mathbf{E}\gamma(x_i, Y_i).$$

Our target function $\theta_n^*$ is defined as

$$(3) \qquad \theta_n^* := \arg\min_{\theta \in \Theta} P\gamma_{f_\theta}.$$

When the target is sparse, i.e., when only a few of the coefficients $\theta_{n,k}^*$ are nonzero, it makes sense to try to prove that also the estimator $\hat\theta_n$ is sparse. Non-asymptotic bounds for this case (albeit with random design) are studied in [12]. It is assumed there that the base functions $\{\psi_k\}$ have design matrix with eigenvalues bounded away from zero (or at least that the base functions corresponding to the non-zero coefficients in $\theta_n^*$ have this property). In the present paper, the base functions are allowed to be highly correlated. We will consider the case where they form a VC class, or more generally, have $\epsilon$-covering number which is polynomial in $1/\epsilon$. This means that a certain smoothness is imposed a priori, and that sparseness is less an issue.

We use the following notation. The empirical distribution based on the sample $\{(x_i, Y_i)\}_{i=1}^n$ is denoted by $P_n$, and the empirical distribution of the covariates $\{x_i\}_{i=1}^n$ is written as $Q_n$. The $L_2(Q_n)$ norm is written as $\|\cdot\|_n$. Moreover, $\|\cdot\|_\infty$ denotes the sup norm (which in our case may be understood as $\|f\|_\infty = \max_{1 \le i \le n} |f(x_i)|$, for a function $f$ on $\mathcal{X}$).

We impose four basic assumptions: Assumptions L, M, A and B.

**Assumption L.** The loss function $\gamma_f$ is of the form $\gamma_f(x,y) = \gamma(f(x), y)$, where $\gamma(\cdot, y)$ is convex for all $y \in \mathcal{Y}$. Moreover, it satisfies the Lipschitz property

$$|\gamma(f_\theta(x), y) - \gamma(f_{\tilde\theta}(x), y)| \le |f_\theta(x) - f_{\tilde\theta}(x)|, \forall\ (x,y) \in \mathcal{X} \times \mathcal{Y},\ \forall\ \theta, \tilde\theta \in \Theta.$$

**Assumption M.** There exists a non-decreasing function $\sigma(\cdot)$, such that all $M > 0$ and all all $\theta \in \Theta$ with $\|f_\theta - f_{\theta_n^*}\|_\infty \le M$, one has

$$P(\gamma_{f_\theta} - \gamma_{f_{\theta_n^*}}) \ge \|f_\theta - f_{\theta_n^*}\|_n^2 / \sigma^2(M).$$

Assumption M thus assumes quadratic margin behavior. In [12], more general margin behavior is allowed, and the choice of the smoothing parameter does not depend on the margin behavior. However, in the setup of the present paper, the choice of the smoothing parameter does depend on the margin behavior.

**Assumption A.** It holds that

$$\|\psi_k\|_\infty \le 1,\ 1 \le k \le m.$$

**Assumption B.** For some constant $A \ge 1$, and for $V = 2/s - 2$, it holds that

$$N(\epsilon, \Psi) \le A\epsilon^{-V}, \forall\ \epsilon > 0.$$

Here $N(\epsilon, \Psi)$ denotes the $\epsilon$-covering number of $(\Psi, \|\cdot\|_n)$, with $\Psi := \{\psi_k\}_{k=1}^m$.

The paper is organized as follows. Section 2 presents a bound for the increments of the empirical process. Section 3 takes such a bound for granted and presents a non-asymptotic bound for $\|f_{\hat\theta_n} - f_{\theta_n^*}\|_n$ and $I(\hat\theta_n)$. The two sections can be read independently. In particular, any improvement of the bound obtained in Section 2 can be directly inserted in the result of Section 3. The proofs, which are perhaps the most interesting part of the paper, are given in Section 4.

## 2. Increments of the empirical process indexed by a subset of a linear space

Let $\varepsilon_1, \ldots, \varepsilon_n$ be i.i.d. random variables, taking values $\pm 1$ each with probability $1/2$. Such a sequence is called a Rademacher sequence. Consider for $\epsilon > 0$ and

$M > 0$, the quantity

$$\mathcal{Z}_{\epsilon,M} := \sup_{\|f_\theta\|_n \leq \epsilon,\ I(\theta) \leq M} \left| \frac{1}{n} \sum_{i=1}^n f_\theta(x_i)\varepsilon_i \right|.$$

We need a bound for the mean $\mathbf{E}\mathcal{Z}_{\epsilon,M}$, because this quantity will occur in the concentration inequality (Theorem 4.1). In [12], the following trivial bound is used:

$$\mathbf{E}\mathcal{Z}_{\epsilon,M} \leq M \mathbf{E}\left( \max_{1 \leq k \leq m} \left| \frac{1}{n} \sum_{i=1}^n \varepsilon_i \psi_k(x_i) \right| \right).$$

On the right hand side, one now has the mean of finitely many functions, which is easily handled (see for example Lemma 4.1). However, when the base functions $\psi_k$ are highly correlated, this bound is too rough. We need therefore to proceed differently.

Let $\mathrm{conv}(\Psi) = \{f_\theta = \sum_{k=1}^m \theta_k \psi_k :\ \theta_k \geq 0, \sum_{k=1}^m \theta_k = 1\}$ be the convex hull of $\Psi$.

Recall that $s = 2/(2+V)$, where $V$ is from Assumption B. From e.g. [10], Lemma 3.2, it can be derived that for some constant $C$, and for all $\epsilon > 0$,

(4)
$$\mathbf{E}\left| \max_{f \in \mathrm{conv}(\Psi), \|f\|_n \leq \epsilon} \frac{1}{n} \sum_{i=1}^n f(x_i)\varepsilon_i \right| \leq C\epsilon^s \frac{1}{\sqrt{n}}.$$

The result follows from the chaining technique, and applying the entropy bound

(5)
$$\log N(\epsilon, \mathrm{conv}(\Psi)) \leq A_0 \epsilon^{-2(1-s)},\ \epsilon > 0,$$

which is derived in [2]. Here, $A_0$ is a constant depending on $V$ and $A$.

**Remark.** It may be verified that the constant $C$ in (4) is then at least proportional to $1/s$, i.e., it is large when $s$ is small.

Our aim is now to obtain a bound from direct calculations. Pollard ([8]) presents the bound

$$\log N(\epsilon, \mathrm{conv}(\Psi)) \leq A_1 \epsilon^{-2(1-s)} \log \frac{1}{\epsilon},\ \epsilon > 0,$$

where $A_1$ is another constant depending on $V$ and $A$. In other words, Pollard's bound has an additional log-factor. On the other hand, we found Pollard's proof a good starting point in our attempt to derive the increments directly, without chaining. This is one of the reasons why our direct bound below has an additional $\log m$ factor. Thus, our result should primarily be seen as illustration that direct calculations are possible.

**Theorem 2.1.** *For $\epsilon \geq 16/m$, and $m \geq 4$, we have*

$$\mathbf{E}\left| \max_{f \in \mathrm{conv}(\Psi), \|f\|_n \leq \epsilon} \frac{1}{n} \sum_{i=1}^n f(x_i)\varepsilon_i \right| \leq 20\sqrt{1+2A}\epsilon^s \sqrt{\frac{\log(6m)}{n}}.$$

Clearly the set $\{\sum_{k=1}^m \theta_k \psi_k :\ I(\theta) \leq 1\}$ is the convex hull of $\{\pm \psi_k\}_{k=1}^m$. Using a renormalization argument, one arrives at the following corollary

**Corollary 2.1.** *We have for $\epsilon/M > 8/m$ and $m \geq 2$*

$$\mathbf{E}\mathcal{Z}_{\epsilon,M} \leq 20\sqrt{1+4A} M^{1-s} \epsilon^s \sqrt{\frac{\log(12m)}{n}}.$$

Invoking symmetrization, contraction and concentration inequalities (see Section 4), we establish the following lemma. We present it in a form convenient for application in the proof of Theorem 3.1.

**Lemma 2.1.** *Define for $\epsilon > 0$, $M > 0$, and $\epsilon/M > 8/m$, $m \geq 2$,*

$$\mathbf{Z}_{\epsilon,M} := \sup_{\|f_\theta - f_{\theta_n^*}\|_n \leq \epsilon,\ I(\theta - \theta_n^*) \leq M} |(P_n - P)(\gamma_{f_\theta} - \gamma_{f_{\theta_n^*}})|.$$

Let

$$\lambda_{n,0} := 80\sqrt{1+4A}\sqrt{\frac{\log(12m)}{n}}.$$

*Then it holds for all $\sigma > 0$, that*

$$\mathbf{P}\left(\mathbf{Z}_{\epsilon,M} \geq \lambda_{n,0}\epsilon^s M^{1-s} + \frac{\epsilon^2}{27\sigma^2}\right) \leq \exp\left[-\frac{n\epsilon^2}{2 \times (27\sigma^2)^2}\right].$$

## 3. A non-asymptotic bound for the estimator

The following theorem presents bounds along the lines of results in [10], [11] and [3], but it is stated in a non-asymptotic form. It moreover formulates explicitly the dependence on the expected increments of the empirical process.

**Theorem 3.1.** *Define for $\epsilon > 0$ and $M > 0$,*

$$\mathbf{Z}_{\epsilon,M} := \sup_{\|f_\theta - f_{\theta_n^*}\|_n \leq \epsilon,\ I(\theta - \theta_n^*) \leq M} |(P_n - P)(\gamma_{f_\theta} - \gamma_{f_{\theta_n^*}})|.$$

*Let $\lambda_{n,0}$ be such that for all $8/m \leq \epsilon/M \leq 1$, we have*

(6) $$\mathbf{E}\mathbf{Z}_{\epsilon,M} \leq \lambda_{n,0}\epsilon^s M^{1-s}.$$

*Let $c \geq 3$ be some constant.*
*Define*

$$M_n := 2^{\frac{2-s}{2(1-s)}}(27)^{-\frac{s}{2(1-s)}} c^{\frac{1}{1-s}} I(\theta_n^*),$$

$$\sigma_n^2 := \sigma^2(M_n),$$

*and*

$$\epsilon_n := \sqrt{54}\sigma_n^{\frac{2}{2-s}} c^{\frac{1}{2-s}} \lambda_{n,0}^{\frac{1}{2-s}} I^{\frac{1-s}{2-s}}(\theta_n^*) \vee 27\sigma_n^2 \lambda_{n,0}.$$

*Assume that*

(7) $$1 \leq \left(\frac{27}{2}\right)^{-\frac{2-s}{2(1-s)}} c^{\frac{1}{1-s}} \frac{1}{\sigma_n^2 \lambda_{n,0}} I(\theta_n^*) \leq \left(\frac{m}{8}\right)^{2-s}.$$

*Then for $\lambda_n := c\sigma_n^s \lambda_{n,0}$, with probability at least*

$$1 - \exp\left[-\frac{n\lambda_{n,0}^{\frac{2}{2-s}} c^{\frac{2}{2-s}} I^{\frac{2(1-s)}{2-s}}(\theta_n^*)}{27\sigma_n^{\frac{4(1-s)}{2-s}}}\right],$$

*we have that*

$$\|f_{\hat{\theta}_n} - f_{\theta_n^*}\|_n \leq \epsilon_n$$

*and*

$$I(\hat{\theta}_n - \theta_n^*) \leq M_n.$$

Let us formulate the asymptotic implication of Theorem 3.1 in a corollary. For positive sequences $\{a_n\}$ and $\{b_n\}$, we use the notation

$$a_n \asymp b_n,$$

when

$$0 < \liminf_{n \to \infty} \frac{a_n}{b_n} \leq \limsup_{n \to \infty} \frac{a_n}{b_n} < \infty.$$

The corollary yields e.g. the rate $\epsilon_n \asymp n^{-1/3}$ for the case where the penalty represents the total variation of a function $f$ on $\{x_1, \ldots, x_n\} \subset \mathbf{R}$ (in which case $s = 1/2$).

**Corollary 3.1.** *Suppose that $A$ and $s$ do not depend on $n$, and that $I(\theta_n^*) \asymp 1$ and $\sigma^2(M_n) \asymp 1$ for all $M_n \asymp 1$. By (4), we may take $\lambda_n \asymp 1/\sqrt{n}$, in which case, with probability $1 - \exp[-d_n]$, it holds that $\|f_{\hat{\theta}_n} - f_{\theta_n^*}\|_n \leq \epsilon_n$, and $I(\hat{\theta}_n - \theta_n^*) \leq M_n$, with*

$$\epsilon_n \asymp n^{-\frac{1}{2(2-s)}}, \ M_n \asymp 1, \ d_n \asymp n\epsilon_n^2 \asymp n^{\frac{1-s}{2-s}}.$$

## 4. Proofs

### 4.1. Preliminaries

**Theorem 4.1** (Concentration theorem [6]). *Let $Z_1, \ldots, Z_n$ be independent random variables with values in some space $\mathcal{Z}$ and let $\Gamma$ be a class of real-valued functions on $\mathcal{Z}$, satisfying*

$$a_{i,\gamma} \leq \gamma(Z_i) \leq b_{i,\gamma},$$

*for some real numbers $a_{i,\gamma}$ and $b_{i,\gamma}$ and for all $1 \leq i \leq n$ and $\gamma \in \Gamma$. Define*

$$L^2 := \sup_{\gamma \in \Gamma} \sum_{i=1}^n (b_{i,\gamma} - a_{i,\gamma})^2/n,$$

*and*

$$\mathbf{Z} := \sup_{\gamma \in \Gamma} \left| \frac{1}{n} \sum_{i=1}^n (\gamma(Z_i) - E\gamma(Z_i)) \right|.$$

*Then for any positive $z$,*

$$\mathbf{P}(\mathbf{Z} \geq \mathbf{EZ} + z) \leq \exp\left[-\frac{nz^2}{2L^2}\right].$$

The Concentration theorem involves the expectation of the supremum of the empirical process. We derive bounds for it using symmetrization and contraction. Let us recall these techniques here.

**Theorem 4.2** (Symmetrization theorem [13]). *Let $Z_1, \ldots, Z_n$ be independent random variables with values in $\mathcal{Z}$, and let $\varepsilon_1, \ldots, \varepsilon_n$ be a Rademacher sequence independent of $Z_1, \ldots, Z_n$. Let $\Gamma$ be a class of real-valued functions on $\mathcal{Z}$. Then*

$$\mathbf{E}\left(\sup_{\gamma \in \Gamma} \left| \sum_{i=1}^n \{\gamma(Z_i) - E\gamma(Z_i)\} \right| \right) \leq 2\mathbf{E}\left(\sup_{\gamma \in \Gamma} \left| \sum_{i=1}^n \varepsilon_i \gamma(Z_i) \right| \right).$$

**Theorem 4.3** (Contraction theorem [5])**.** *Let $z_1, \ldots, z_n$ be non-random elements of some space $\mathcal{Z}$ and let $\mathcal{F}$ be a class of real-valued functions on $\mathcal{Z}$. Consider Lipschitz functions $\gamma_i : \mathbf{R} \to \mathbf{R}$, i.e.*

$$|\gamma_i(s) - \gamma_i(\tilde{s})| \leq |s - \tilde{s}|, \ \forall \ s, \tilde{s} \in \mathbf{R}.$$

*Let $\varepsilon_1, \ldots, \varepsilon_n$ be a Rademacher sequence. Then for any function $f^* : \mathcal{Z} \to \mathbf{R}$, we have*

$$\mathbf{E}\left(\sup_{f \in \mathcal{F}} \left|\sum_{i=1}^n \varepsilon_i \{\gamma_i(f(z_i)) - \gamma_i(f^*(z_i))\}\right|\right) \leq 2\mathbf{E}\left(\sup_{f \in \mathcal{F}} \left|\sum_{i=1}^n \varepsilon_i(f(z_i) - f^*(z_i))\right|\right).$$

We now consider the case where $\Gamma$ is a finite set of functions.

**Lemma 4.1.** *Let $Z_1, \ldots, Z_n$ be independent $\mathcal{Z}$-valued random variables, and $\gamma_1, \ldots, \gamma_m$ be real-valued functions on $\mathcal{Z}$, satisfying*

$$a_{i,k} \leq \gamma_k(Z_i) \leq b_{i,k},$$

*for some real numbers $a_{i,k}$ and $b_{i,k}$ and for all $1 \leq i \leq n$ and $1 \leq k \leq m$. Define*

$$L^2 := \max_{1 \leq k \leq m} \sum_{i=1}^n (b_{i,k} - a_{i,k})^2 / n,$$

*Then*

$$\mathbf{E}\left(\max_{1 \leq k \leq m} \left|\frac{1}{n} \sum_{i=1}^n \{\gamma_k(Z_i) - E\gamma_k(Z_i)\}\right|\right) \leq 2L\sqrt{\frac{\log(3m)}{n}}.$$

*Proof.* The proof uses standard arguments, as treated in e.g. [13]. Let us write for $1 \leq k \leq m$,

$$\bar{\gamma}_k := \frac{1}{n} \sum_{i=1}^n \left\{\gamma_k(Z_i) - E\gamma_k(Z_i)\right\}.$$

By Hoeffding's inequality, for all $z \geq 0$

$$\mathbf{P}(|\bar{\gamma}_k| \geq z) \leq 2 \exp\left[-\frac{nz^2}{2L^2}\right].$$

Hence,

$$\mathbf{E}\exp\left[\frac{n}{4L^2}\bar{\gamma}_k^2\right] = 1 + \int_1^\infty \mathbf{P}\left(|\bar{\gamma}_k| \geq \sqrt{\frac{4L^2}{n}\log t}\right) dt$$

$$\leq 1 + 2\int_1^\infty \frac{1}{t^2} dt = 3.$$

Thus

$$\mathbf{E}\left(\max_{1 \leq k \leq m} |\bar{\gamma}_k|\right) = \frac{2L}{\sqrt{n}} \mathbf{E}\sqrt{\max_{1 \leq k \leq m} \log \exp\left[\frac{4}{4L^2}\bar{\gamma}_k^2\right]}$$

$$\leq \frac{2L}{\sqrt{n}}\sqrt{\log \mathbf{E}\left(\max_{1 \leq k \leq m} \exp\left[\frac{4}{4L^2}\bar{\gamma}_k^2\right]\right)} \leq 2L\sqrt{\frac{\log(3m)}{n}}. \quad \square$$

## 4.2. Proofs of the results in Section 2

*Proof of Theorem 2.1.* Let us define, for $k = 1, \ldots, m$,

$$\xi_k := \frac{1}{n} \sum_{i=1}^{n} \psi_k(x_i) \varepsilon_i.$$

We have

$$\frac{1}{n} \sum_{i=1}^{n} f_\theta(x_i) \varepsilon_i = \sum_{k=1}^{m} \theta_k \xi_k.$$

Partition $\{1, \ldots, m\}$ into $N := N(\epsilon^s, \Psi)$ sets $V_j$, $j = 1, \ldots, N$, such that

$$\|\psi_k - \psi_l\|_n \leq 2\epsilon^s, \ \forall \ k, l \in V_j.$$

We can write

$$\frac{1}{n} \sum_{i=1}^{n} f_\theta(x_i) \varepsilon_i = \sum_{j=1}^{N} \alpha_j \sum_{k \in V_j} p_{j,k} \xi_k,$$

where

$$\alpha_j = \alpha_j(\theta) := \sum_{k \in V_j} \theta_k, \ p_{j,k} = p_{j,k}(\theta) := \frac{\theta_k}{\alpha_j}.$$

Set for $j = 1, \ldots, N$,

$$n_j = n_j(\alpha) := 1 + \lfloor \frac{\alpha_j}{\epsilon^{2(1-s)}} \rfloor.$$

Choose $\pi_{t,j} = \pi_{t,j}(\theta)$, $t = 1, \ldots, n_j$, $j = 1, \ldots, N$ independent random variables, independent of $\varepsilon_1, \ldots, \varepsilon_n$, with distribution

$$\mathbf{P}(\pi_{t,j} = k) = p_{j,k}, \ k \in V_j, \ j = 1, \ldots, N.$$

Let $\bar{\psi}_j = \bar{\psi}_j(\theta) := \sum_{i=1}^{n_j} \psi_{\pi_{t,j}}/n_j$ and $\bar{\xi}_j = \bar{\xi}_j(\theta) := \sum_{i=1}^{n_j} \xi_{\pi_{t,j}}/n_j$.

We will choose a realization $\{(\psi_j^*, \xi_j^*) = (\psi_j^*(\theta), \xi_j^*(\theta))\}_{j=1}^N$ of $\{(\bar{\psi}_j, \bar{\xi}_j)\}_{j=1}^N$ depending on $\{\varepsilon_i\}_{i=1}^n$, satisfying appropriate conditions (namely, (9) and (10) below). We may then write

$$\left| \sum_{k=1}^{m} \theta_k \xi_k \right| \leq \left| \sum_{j=1}^{N} \alpha_j \xi_j^* \right| + \left| \sum_{k=1}^{m} \theta_k \xi_k - \sum_{j=1}^{N} \alpha_j \xi_j^* \right|.$$

Consider now

$$\sum_{j=1}^{N} \alpha_j \xi_j^*.$$

Let $\mathcal{A}^N := \{\sum_{i=1}^{N} \alpha_j = 1, \ \alpha_j \geq 0\}$. Endow $\mathcal{A}^N$ with the $\ell_1$ metric. The $\epsilon$-covering number $D(\epsilon)$ of $\mathcal{A}^N$ satisfies the bound

$$D(\epsilon) \leq \left(\frac{4}{\epsilon}\right)^N.$$

Let $\mathcal{A}_\epsilon$ be a maximal $\epsilon$-covering set of $\mathcal{A}^N$. For all $\alpha \in \mathcal{A}$ there is an $\alpha' \in \mathcal{A}_\epsilon$ such that $\sum_{j=1}^{N} |\alpha_j - \alpha'_j| \leq \epsilon$.

We now write

$$\left|\sum_{k=1}^{m}\theta_k\xi_k\right| \leq \left|\sum_{j=1}^{N}(\alpha_j - \alpha'_j)\xi_j^*\right| + \left|\sum_{k=1}^{m}\theta_k\xi_k - \sum_{j=1}^{N}\alpha_j\xi_j^*\right| + \left|\sum_{j=1}^{N}\alpha'_j\xi_j^*\right|$$
$$:= \mathrm{i}(\theta) + \mathrm{ii}(\theta) + \mathrm{iii}(\theta).$$

Let $\Pi$ be the set of possible values of the vector $\{\pi_{t,j} : t = 1,\ldots,n_j, j = 1,\ldots,N\}$, as $\theta$ varies. Clearly,

$$\mathrm{i}(\theta) \leq \epsilon \max_{\Pi} \max_{j} |\bar{\xi}_j|,$$

where we take the maximum over all possible realizations of $\{\bar{\xi}_j\}_{j=1}^{N}$ over all $\theta$.

For each $t$ and $j$, $\pi_{t,j}$ takes its values in $\{1,\ldots,m\}$, that is, it takes at most $m$ values. We have

$$\sum_{j=1}^{N} n_j \leq N + \sum_{j=1}^{N} \frac{\alpha_j}{\epsilon^{2(1-s)}}$$
$$\leq A\epsilon^{-sV} + \frac{\sum_{k=1}^{m}\theta_k}{\epsilon^{2(1-s)}}$$
$$= (1+A)\epsilon^{-2(1-s)} \leq K+1.$$

where $K$ is the integer

$$K := \lfloor (1+A)\epsilon^{2(1-s)} \rfloor.$$

The number of integer sequences $\{n_j\}_{j=1}^{N}$ with $\sum_{j=1}^{N} n_j \leq K+1$ is equal to

$$\binom{N+K+2}{K+1} \leq 2^{N+K+2} \leq 4 \times 2^{(1+2A)\epsilon^{-2(1-s)}}.$$

So the cardinality $|\Pi|$ of $\Pi$ satisfies

$$|\Pi| \leq 4 \times 2^{(1+2A)\epsilon^{-2(1-s)}} \times m^{(1+A)\epsilon^{-2(1-s)}} \leq (2m)^{(1+2A)\epsilon^{-2(1-s)}},$$

since $A \geq 1$ and $m \geq 4$.

Now, since $\|\psi\|_\infty \leq 1$ for all $\psi \in \Psi$, we know that for any convex combination $\sum_k p_k\xi_k$, one has $\mathbf{E}|\sum_k p_k\xi_k|^2 \leq 1/n$. Hence $\mathbf{E}\bar{\xi}_j^2 \leq 1/n$ for any fixed $\bar{\xi}_j$ and thus, by Lemma 4.1,

(8) $\quad \epsilon\mathbf{E}\max_{\Pi}\max_{j}|\bar{\xi}_j| \leq 2\epsilon\sqrt{1+2A}\epsilon^{-(1-s)}\sqrt{\frac{\log(6m)}{n}} = 2\sqrt{1+2A}\epsilon^s\sqrt{\frac{\log(6m)}{n}}.$

We now turn to $\mathrm{ii}(\theta)$.

By construction, for $i = 1,\ldots,n$, $t = 1,\ldots,n_j$, $j = 1,\ldots,N$,

$$\mathbf{E}\psi_{\pi_{t,j}}(x_i) = \sum_{k \in V_j} p_{j,k}\psi_k(x_i) := g_j(x_i)$$

and hence

$$\mathbf{E}(\psi_{\pi_{t,j}}(x_i) - g_j(x_i))^2 \leq \max_{k,l \in V_j} (\psi_k(x_i) - \psi_l(x_i))^2.$$

Thus

$$\mathbf{E}(\bar{\psi}_j(x_i) - g_j(x_i))^2 \leq \max_{k,l \in V_j} (\psi_k(x_i) - \psi_l(x_i))^2/n_j,$$

and so
$$\mathbf{E}\|\bar{\psi}_j - g_j\|_n^2 \leq \max_{k,l \in V_j} \|\psi_k - \psi_l\|_n^2/n_j \leq (2\epsilon^s)^2/n_j = 4\epsilon^{2s}/n_j.$$

Therefore
$$\mathbf{E}\left\|\sum_{j=1}^N \alpha_j(\bar{\psi}_j - g_j)\right\|_n^2 = \sum_{j=1}^N \alpha_j^2 \mathbf{E}\|\bar{\psi}_j - g_j\|_n^2$$
$$\leq 4\epsilon^{2s} \sum_{j=1}^N \frac{\alpha_j^2}{n_j} \leq 4\epsilon^{2s} \sum_{j=1}^N \frac{\alpha_j^2 \epsilon^{2(1-s)}}{\alpha_j} \leq 4\epsilon^2.$$

Let $\mathbf{E}_\varepsilon$ denote conditional expectation given $\{\varepsilon_i\}_{i=1}^n$. Again, by construction
$$\mathbf{E}_\varepsilon \xi_{\pi_{t,j}} = \sum_{k \in V_j} p_{j,k} \xi_k := e_j = e_j(\theta),$$

and hence
$$\mathbf{E}_\varepsilon(\xi_{\pi_{t,j}} - e_j)^2 \leq \max_{k,l \in V_j} (\xi_k - \xi_l)^2.$$

Thus
$$\mathbf{E}_\varepsilon(\bar{\xi}_j - e_j)^2 \leq \max_{k,l \in V_j} (\xi_k - \xi_l)^2/n_j.$$

So we obtain
$$\mathbf{E}_\varepsilon\left|\sum_{j=1}^N \alpha_j(\bar{\xi}_j - e_j)\right| \leq \sum_{j=1}^N \alpha_j \mathbf{E}_\varepsilon|\bar{\xi}_j - e_j| \leq \sum_{j=1}^N \alpha_j \max_{k,l \in V_j} \frac{|\xi_k - \xi_l|}{\sqrt{n_j}}$$
$$\leq \sum_{j=1}^N \frac{\alpha_j \epsilon^{1-s}}{\sqrt{\alpha_j}} \max_{k,l \in V_j} |\xi_k - \xi_l| = \sum_{j=1}^N \sqrt{\alpha_j} \epsilon^{1-s} \max_{k,l \in V_j} |\xi_k - \xi_l|$$
$$\leq \sqrt{N} \epsilon^{1-s} \max_j \max_{k,l \in V_j} |\xi_k - \xi_l| \leq \sqrt{A} \max_j \max_{k,l \in V_j} |\xi_k - \xi_l|.$$

It follows that, given $\{\varepsilon_i\}_{i=1}^n$, there exists a realization
$$\{(\psi_j^*, \xi_j^*) = (\psi_j^*(\theta), \xi_j^*(\theta))\}_{j=1}^N$$
of $\{(\bar{\psi}_j, \bar{\xi}_j)\}_{j=1}^N$ such that

(9) $$\|\sum_{j=1}^N \alpha_j(\psi_j^* - g_j)\|_n^2 \leq 4\epsilon$$

as well as

(10) $$\left|\sum_{j=1}^N \alpha_j(\xi_j^* - e_j)\right| \leq 2\sqrt{A} \max_j \max_{k,l \in V_j} |\xi_k - \xi_l|.$$

Thus we have
$$\text{ii}(\theta) \leq 2\sqrt{A} \max_j \max_{k,l \in V_j} |\xi_k - \xi_l|.$$

Since $\mathbf{E}|\xi_k - \xi_l|^2 \leq 2\epsilon^2/n$ for all $k, l \in V_j$ and all $j$, we have by Lemma 4.1,

(11) $$2\sqrt{A}\mathbf{E} \max_j \max_{k,l \in V_j} |\xi_k - \xi_l| \leq 6\sqrt{A}\epsilon^s \sqrt{\frac{\log(6m)}{n}}.$$

Finally, consider iii($\theta$). We know that

$$\|f_\theta\|_n = \left\|\sum_{j=1}^{N} \alpha_j g_j\right\|_n \leq \epsilon.$$

Moreover, we have shown in (9) that

$$\left\|\sum_{j=1}^{N} \alpha_j(\psi_j^* - g_j)\right\|_n \leq 4\epsilon.$$

Also

$$\left\|\sum_{j=1}^{N} (\alpha_j - \alpha_j')\psi_j^*\right\|_n \leq \sum_{j=1}^{N} |\alpha_j - \alpha_j'|\|\psi_j^*\|_n \leq \epsilon,$$

since $\|\psi_j^*\|_\infty \leq 1$ for all $j$. Thus

$$\left\|\sum_{j=1}^{N} \alpha_j'\psi_j^*\right\|_n \leq \left\|\sum_{j=1}^{N}(\alpha_j' - \alpha_j)\psi_j^*\right\|_n + \left\|\sum_{j=1}^{N}\alpha_j(\psi_j^* - g_j)\right\|_n + \|f_\theta\|_n \leq 6\epsilon.$$

The total number of functions of the form $\sum_{j=1}^{N} \alpha_j'\xi_j^*$ is bounded by

$$\left(\frac{4}{\epsilon}\right)^N \times |\Pi| \leq \left(\frac{4}{\epsilon}\right)^{A\epsilon^{-2(1-s)}} \times (2m)^{(1+2A)\epsilon^{-2(1-s)}}$$
$$\leq (2m)^{(1+2A)\epsilon^{-2(1-s)}},$$

since we assume $\epsilon \geq 16/m$, and $A \geq 1$. Hence, by Lemma 4.1,

(12) $$\mathbf{E}\max_{\alpha' \in \mathcal{A}_\epsilon}\max_{\Pi}|\sum_{j=1}^{N}\alpha_j'\xi_j^*| \leq 12\sqrt{1+2A}\epsilon^s\sqrt{\frac{\log(6m)}{n}}.$$

We conclude from (8), (11), and (12), that

$$\mathbf{E}\max_{\theta}\left|\sum_{j=1}^{N}\alpha_j(\theta)e_j(\theta)\right|$$
$$\leq 2\sqrt{1+2A}\epsilon^s\sqrt{\frac{\log(6m)}{n}} + 6\sqrt{A}\epsilon^s\sqrt{\frac{\log(6m)}{n}} + 12\sqrt{1+2A}\epsilon^s\sqrt{\frac{\log(6m)}{n}}$$
$$\leq 20\sqrt{1+2A}\epsilon^s\sqrt{\frac{\log(6m)}{n}}. \qquad \square$$

*Proof of Lemma 2.1.* Let

$$\mathcal{Z}_{\epsilon,M} := \sup_{\|f_\theta\|_n \leq \epsilon,\, I(\theta) \leq M}\left|\frac{1}{n}\sum_{i=1}^{n}\gamma_{f_\theta}(x_i, Y_i)\varepsilon_i\right|$$

denote the symmetrized process. Clearly, $\{f_\theta = \sum_{k=1}^{m}\theta_k\psi_k :\ I(\theta) = 1\}$ is the convex hull of $\tilde{\Psi} := \{\pm\psi_k\}_{k=1}^{m}$. Moreover, we have

$$N(\epsilon, \tilde{\Psi}) \leq 2N(\epsilon, \Psi).$$

Now, apply Theorem 2.1, to $\tilde{\Psi}$, and use a rescaling argument, to see that

$$\mathbf{E}\mathcal{Z}_{\epsilon,M} \leq 20\sqrt{1+4A}\epsilon^s M^{1-s}\sqrt{\frac{\log(12m)}{n}}.$$

Then from Theorem 4.2 and Theorem 4.3, we know that

$$\mathbf{E}Z_{\epsilon,M} \leq 4\mathbf{E}\mathcal{Z}_{\epsilon,M}.$$

The result now follows by applying Theorem 4.1. □

### 4.3. Proofs of the results in Section 3

The proof of Theorem 3.1 depends on the following simple convexity trick.

**Lemma 4.2.** Let $\epsilon > 0$ and $M > 0$. Define $\tilde{f}_n = t\hat{f}_n + (1-t)f_n^*$ with

$$t := (1 + \|\hat{f}_n - f_n^*\|_n/\epsilon + I(\hat{f}_n - f_n^*)/M)^{-1},$$

and with $\hat{f}_n := f_{\hat{\theta}_n}$ and $f_n^* := f_{\theta_n^*}$. When it holds that

$$\|\tilde{f}_n - f_n^*\|_n \leq \frac{\epsilon}{3}, \text{ and } I(\tilde{f}_n - f_n^*) \leq \frac{M}{3},$$

then

$$\|\hat{f}_n - f_n^*\|_n \leq \epsilon, \text{ and } I(\hat{f}_n - f_n^*) \leq M.$$

*Proof.* We have

$$\tilde{f}_n - f_n^* = t(\hat{f}_n - f_n^*),$$

so $\|\tilde{f}_n - f_n^*\|_n \leq \epsilon/3$ implies

$$\|\hat{f}_n - f_n^*\|_n \leq \frac{\epsilon}{3t} = (1 + \|\hat{f}_n - f_n^*\|_n/\epsilon + I(\hat{f}_n - f_n^*)/M)\frac{\epsilon}{3}.$$

So then

(13) $$\|\hat{f}_n - f_n^*\|_n \leq \frac{\epsilon}{2} + \frac{\epsilon}{2M}I(\hat{f}_n - f_n^*).$$

Similarly, $I(\tilde{f}_n - f_n^*) \leq M/3$ implies

(14) $$I(\hat{f}_n - f_n^*) \leq \frac{M}{2} + \frac{M}{2\epsilon}\|\hat{f}_n - f_n^*\|_n.$$

Inserting (14) into (13) gives

$$\|\hat{f}_n - f_n^*\|_n \leq \frac{3\epsilon}{4} + \frac{1}{4}\|\hat{f}_n - f_n^*\|_n,$$

i.e., $\|\hat{f}_n - f_n^*\|_n \leq \epsilon$. Similarly, Inserting (13) into (14) gives $I(\hat{f}_n - f_n^*) \leq M$. □

*Proof of Theorem 3.1.* Note first that, by the definition of of $M_n$, $\epsilon_n$ and $\lambda_n$, it holds that

(15) $$\lambda_{n,0}\epsilon_n^s M_n^{1-s} = \frac{\epsilon_n^2}{27\sigma_n^2},$$

and also

(16) $$(27)^{\frac{s}{2-s}} c^{-\frac{2}{2-s}} \lambda_n^{\frac{2}{2-s}} M_n^{\frac{2(1-s)}{2-s}} = \frac{\epsilon_n^2}{27\sigma_n^2}.$$

Define
$$\tilde{\theta}_n = t\hat{\theta}_n + (1-t)\theta_n^*,$$
where
$$t := (1 + \|f_{\hat{\theta}_n} - f_{\theta_n^*}\|_n/\epsilon_n + I(f_{\hat{\theta}_n} - f_{\theta_n^*})/M_n)^{-1}.$$

We know that by convexity, and since $\hat{\theta}_n$ minimizes the penalized empirical risk, we have

$$P_n \gamma_{f_{\tilde{\theta}_n}} + \lambda_n^{\frac{2}{2-s}} I^{\frac{2(1-s)}{2-s}}(\tilde{\theta}_n)$$
$$\leq t\left(P_n \gamma_{f_{\hat{\theta}_n}} + \lambda_n^{\frac{2}{2-s}} I^{\frac{2(1-s)}{2-s}}(\hat{\theta}_n)\right) + (1-t)\left(P_n \gamma_{f_{\theta_n^*}} + \lambda_n^{\frac{2}{2-s}} I^{\frac{2(1-s)}{2-s}}(\theta_n^*)\right)$$
$$\leq P_n \gamma_{f_{\theta_n^*}} + \lambda_n^{\frac{2}{2-s}} I^{\frac{2(1-s)}{2-s}}(\theta_n^*).$$

This can be rewritten as

$$P(\gamma_{f_{\tilde{\theta}_n}} - \gamma_{f_{\theta_n^*}}) + \lambda_n^{\frac{2}{2-s}} I^{\frac{2(1-s)}{2-s}}(\tilde{\theta}_n) \leq -(P_n - P)(\gamma_{f_{\tilde{\theta}_n}} - \gamma_{f_{\theta_n^*}}) + \lambda_n^{\frac{2}{2-s}} I^{\frac{2(1-s)}{2-s}}(\theta_n^*).$$

Since $I(f_{\tilde{\theta}_n} - f_{\theta_n^*}) \leq M_n$, and $\|\psi_k\|_\infty \leq 1$ (by Assumption A), we have that $\|f_{\tilde{\theta}_n} - f_{\theta_n^*}\|_\infty \leq M_n$. Hence, by Assumption M,

$$P(\gamma_{f_{\tilde{\theta}_n}} - \gamma_{f_{\theta_n^*}}) \geq \|f_{\tilde{\theta}_n} - f_{\theta_n^*}\|_n^2/\sigma_n^2.$$

We thus obtain

$$\frac{\|f_{\tilde{\theta}_n} - f_{\theta_n^*}\|_n^2}{\sigma_n^2} + \lambda_n^{\frac{2}{2-s}} I^{\frac{2(1-s)}{2-s}}(\tilde{\theta}_n - \theta_n^*)$$
$$\leq \frac{\|f_{\tilde{\theta}_n} - f_{\theta_n^*}\|_n^2}{\sigma_n^2} + \lambda_n^{\frac{2}{2-s}} I^{\frac{2(1-s)}{2-s}}(\tilde{\theta}_n) + \lambda_n^{\frac{2}{2-s}} I^{\frac{2(1-s)}{2-s}}(\theta_n^*)$$
$$\leq -(P_n - P)(\gamma_{f_{\tilde{\theta}_n}} - \gamma_{f_{\theta_n^*}}) + 2\lambda_n^{\frac{2}{2-s}} I^{\frac{2(1-s)}{2-s}}(\theta_n^*).$$

Now, $\|f_{\tilde{\theta}_n} - f_{\theta_n^*}\|_n \leq \epsilon_n$ and $I(\tilde{\theta}_n - \theta_n^*) \leq M_n$. Moreover $\epsilon_n/M_n \leq 1$ and in view of (7), $\epsilon_n/M_n \geq 8/m$. Therefore, we have by (6) and Theorem 4.1, with probability at least
$$1 - \exp\left[-\frac{n\epsilon_n^2}{2 \times (27\sigma_n^2)^2}\right],$$
that

$$\frac{\|f_{\tilde{\theta}_n} - f_{\theta_n^*}\|_n^2}{\sigma_n^2} + \lambda_n^{\frac{2}{2-s}} I^{\frac{2(1-s)}{2-s}}(\tilde{\theta}_n - \theta_n^*)$$
$$\leq \lambda_{n,0} \epsilon_n^s M_n^{1-s} + 2\lambda_n^{\frac{2}{2-s}} I^{\frac{2(1-s)}{2-s}}(\theta_n^*) + \frac{\epsilon_n^2}{27\sigma_n^2}$$
$$\leq \lambda_{n,0} \epsilon_n^s M_n^{1-s} + (27)^{\frac{s}{2-s}} c^{-\frac{2}{2-s}} \lambda_n^{\frac{2}{2-s}} M_n^{\frac{2(1-s)}{2-s}} + \frac{\epsilon_n^2}{27\sigma_n^2}$$
$$= \frac{1}{9\sigma_n^2} \epsilon_n^2,$$

where in the last step, we invoked (15) and (16).

It follows that
$$\|f_{\tilde{\theta}_n} - f_{\theta_n^*}\|_n \leq \frac{\epsilon_n}{3},$$

and also that
$$I^{\frac{2(1-s)}{2-s}}(\tilde{\theta}_n - \theta_n^*) \leq \frac{\epsilon_n^2}{9\sigma_n^2}\lambda_n^{-\frac{2}{2-s}} \leq \left(\frac{M_n}{3}\right)^{\frac{2(1-s)}{2-s}},$$

since $c \geq 3$.

To conclude the proof, apply Lemma 4.2. □

## References

[1] ALEXANDER, K. S. (1985). Rates of growth for weighted empirical processes. *Proc. Berkeley Conf. in Honor of Jerzy Neyman and Jack Kiefer* **2** 475–493. University of California Press, Berkeley. MR0822047

[2] BALL, K. AND PAJOR, A. (1990). The entropy of convex bodies with "few" extreme points. *Geometry of Banach Spaces (Strobl., 1989)* 25–32. London Math. Soc. Lecture Note Ser. **158**. Cambridge Univ. Press. MR1110183

[3] BLANCHARD, G., LUGOSI, G. AND VAYATIS, N. (2003). On the rate of convergence of regularized boosting classifiers. *J. Machine L. Research* **4** 861–894. MR2076000

[4] BOUSQUET, O. (2002). A Bennett concentration inequality and its application to suprema of empirical processes. *C. R. Acad. Sci. Paris* **334** 495–500. MR1890640

[5] LEDOUX, M. and TALAGRAND, M. (1991). *Probability in Banach Spaces: Isoperimetry and Processes*. Springer, New York. MR1102015

[6] MASSART, P. (2000). About the constants in Talagrand's concentration inequalities for empirical processes. *Ann. Probab.* **28** 863–884. MR1782276

[7] MASSART, P. (2000). Some applications of concentration inequalities to statistics. *Ann. Fac. Sci. Toulouse* **9** 245–303. MR1813803

[8] POLLARD, D. (1990). *Empirical Processes: Theory and Applications*. IMS, Hayward, CA. MR1089429

[9] TALAGRAND, M. (1995). Concentration of measure and isoperimetric inequalities in product spaces. *Publ. Math. de l'I.H.E.S.* **81** 73–205. MR1361756

[10] VAN DE GEER, S. (2000). *Empirical Processes in M-Estimation*. Cambridge Univ. Press.

[11] VAN DE GEER, S. (2002). M-estimation using penalties or sieves. *J. Statist. Planning Inf.* **108** 55–69. MR1947391

[12] VAN DE GEER, S. (2006). High-dimensional generalized linear models and the Lasso. Research Report 133, Seminar für Statistik, ETH Zürich *Ann. Statist.* To appear.

[13] VAN DER VAART, A. W. AND WELLNER, J. A. (1996). *Weak Convergence and Empirical Processes*. Springer, New York. MR1385671